\documentclass{elsart}
\usepackage[all]{xy}
\usepackage{ifthen}
\usepackage{algorithm}
\usepackage{algorithmic}

\newcommand{\cyc}{\textup{\bf c}}
\newcommand{\dec}{\textup{\bf d}}
\newcommand{\len}{\textup{len}}
\newcommand{\sinf}{\inf_\textup{\scriptsize s}}
\newcommand{\ssup}{\sup_\textup{\scriptsize s}}
\newcommand{\slen}{\len_\textup{\scriptsize s}}
\newcommand{\Magma}{\textsc{Magma}}
\newcommand{\proof}{{\textit Proof:\hspace{1em}}}
\renewcommand{\phi}{\varphi}

\newtheorem{mytheorem}{Theorem}[section]
\newenvironment{theorem}[1][]
     {\smallskip\ifthenelse{\equal{#1}{}}{\begin{mytheorem}}
                                         {\begin{mytheorem}[#1]}}
     {\end{mytheorem}\smallskip}
\newtheorem{mynotation}[mytheorem]{Notation}
\newenvironment{notation}[1][]
     {\smallskip\ifthenelse{\equal{#1}{}}{\begin{mynotation}}
                                         {\begin{mynotation}[#1]}\em}
     {\end{mynotation}\smallskip}
\newtheorem{mydefinition}[mytheorem]{Definition}
\newenvironment{definition}[1][]
     {\smallskip\ifthenelse{\equal{#1}{}}{\begin{mydefinition}}
                                         {\begin{mydefinition}[#1]}\em}
     {\end{mydefinition}\smallskip}
\newtheorem{myexample}[mytheorem]{Example}
\newenvironment{example}[1][]
     {\smallskip\ifthenelse{\equal{#1}{}}{\begin{myexample}}
                                         {\begin{myexample}[#1]}\em}
     {\end{myexample}\smallskip}
\newtheorem{myremark}[mytheorem]{Remark}
\newenvironment{remark}[1][]
     {\smallskip\ifthenelse{\equal{#1}{}}{\begin{myremark}}
                                         {\begin{myremark}[#1]}\em}
     {\end{myremark}\smallskip}
\newtheorem{mycorollary}[mytheorem]{Corollary}
\newenvironment{corollary}[1][]
     {\smallskip\ifthenelse{\equal{#1}{}}{\begin{mycorollary}}
                                         {\begin{mycorollary}[#1]}}
     {\end{mycorollary}\smallskip}
\newtheorem{mylemma}[mytheorem]{Lemma}
\newenvironment{lemma}[1][]
     {\smallskip\ifthenelse{\equal{#1}{}}{\begin{mylemma}}
                                         {\begin{mylemma}[#1]}}
     {\end{mylemma}\smallskip}
\newtheorem{myproposition}[mytheorem]{Proposition}
\newenvironment{proposition}[1][]
     {\smallskip\ifthenelse{\equal{#1}{}}{\begin{myproposition}}
                                         {\begin{myproposition}[#1]}}
     {\end{myproposition}\smallskip}
\newtheorem{myalgo}[mytheorem]{Algorithm}
\newenvironment{algo}[1][]
     {\smallskip\ifthenelse{\equal{#1}{}}{\begin{myalgo}}
                                         {\begin{myalgo}[#1]}}
     {\end{myalgo}\smallskip}
\newenvironment{oldtheorem}[1]
     {\smallskip{\bf Theorem #1}\em}
     {\smallskip\par}

\begin{document}

\begin{frontmatter}

\title{A New Approach to the Conjugacy Problem in Garside Groups}
\author{Volker Gebhardt}
\address{School of Mathematics and Statistics, University of Sydney, Australia}
\ead{volker@maths.usyd.edu.au}

\begin{abstract}
 The cycling operation
 endows the super summit set
 $S_x$ of any element $x$ of a Garside group $G$ with the structure of a
 directed graph $\Gamma_x$.
 We establish that the subset $U_x$ of $S_x$ consisting of the circuits
 of $\Gamma_x$ can be used instead of $S_x$ for deciding conjugacy to $x$
 in $G$, yielding a faster and more practical
 solution to the conjugacy problem for Garside
 groups. Moreover, we present a probabilistic approach to the conjugacy
 search problem in Garside groups. The results are likely to have implications
 for the security of recently proposed cryptosystems based on the hardness
 of problems related to the conjugacy (search) problem in braid groups.
\end{abstract}

\end{frontmatter}

\section{Introduction}
\label{sect_intro}

Given a group $G$, the {\em conjugacy problem} in $G$ is to decide for given
elements $a,b\in G$, whether $a$ and $b$ are conjugate in $G$, that is,
whether there exists an element $c\in G$ such that $a^c=b$.  The {\em conjugacy
search problem} in $G$, on the other hand, is to find for given elements
$a$ and $b$ which are known to be conjugate in $G$, an element $c\in G$ such
that $a^c=b$.

Both problems are known to be solvable in Garside groups, that is, in
particular in braid groups
\cite{braid_garside_alg,braid_positive_braids,braid_word_problem,braid_conjugacy_gaussian,braid_conjugacy}.  However, all known algorithms involve
computing a particular invariant of the conjugacy class, the so-called
{\em super summit set}, for either $a$ or $b$ and both the memory and the
time complexity of these algorithms are proportional to the cardinality of
this set.  In the case of the braid group $B_n$,
the best proven bound for this cardinality is exponential in both the braid
index $n$ and the element length $r$ and, while the existence of polynomial
bounds is conjectured, computations in practice are hard or infeasible
even for moderate values of $n$ and $r$.

Recently, braid groups came under interest as possible sources for public
key cryptosystems and the security of most of the proposed cryptosystems
depends on the hardness of variations of the conjugacy (search) problem
\cite{algebraic_public_key,braid_cryptosystem}. Hence an improved understanding
of the conjugacy problems is highly desirable.

The crucial point in computing the super summit set $S_x$ of an element $x$
is the following \lq\lq{}convexity\rq\rq{} property:  For any pair of elements
$u,v\in S_x$ there are elements $u_0,\dots,u_k$ with $u_0 = u$ and
$u_k = v$, such that for $i=1,\dots,k$, $u_i$ is obtained from $u_{i-1}$
by conjugation with a suitable element from a finite set $D$.  This allows us
to compute $S_x$, starting with a single representative, as closure with
respect to conjugation by elements of $D$.

In this paper we establish that a subset of the super summit set, which
in general is much smaller, can be used for deciding
conjugacy in Garside groups.  The set $S_x$ can be endowed with the structure
of a directed graph and we will show that the union of the circuits of
this graph has the same \lq\lq{}convexity\rq\rq{} property as described above,
that is, can be computed in a similar way.  The graph structure used for
proving this result also yields a fast probabilistic algorithm for solving
the conjugacy search problem.

\subsection{Garside Groups and Monoids}
\label{sect_garside_groups}

We start with a brief review of some basic terminology and facts about Garside
groups.  The results can be found, for example, in \cite{braid_garside_alg,braid_positive_braids,braid_word_problem,braid_conjugacy_gaussian,braid_epsteinetal,braid_gaussian_garside,braid_garside_groups}.
Throughout this section, let $M$ be a (left and right) cancellative monoid.

\begin{definition}
\label{def_partial_ordering}
We define partial orderings $\preceq$ and $\succeq$
on the elements of $M$ as follows.  For $a,b\in M$ we say $a\preceq b$
if there exists an element $c\in M$ such that $ac=b$ and we say $a\succeq b$
if there exists an element $c\in M$ such that $a=cb$.

We call $m$ a
{\em (left) lcm} of $a$ and $b$ if $a\preceq m$, $b\preceq m$ and if for
any $x\in M$, $a\preceq x$ and $b\preceq x$ implies $m\preceq x$.
Similarly, we call $d$ a {\sl (left) gcd} of $a$ and $b$ if $d\preceq a$,
$d\preceq b$ and if for any $x\in M$, $x\preceq a$ and $x\preceq b$
implies $x\preceq d$.

\end{definition}

\begin{definition}
\label{def_atom}
$x\in M$ is called an {\em atom} if $x\ne 1$ and if
$x=ab$ for $a,b\in M$ implies $a=1$ or $b=1$.  $M$ is called {\em atomic}
if $M$ is generated by its atoms and if for every $a\in M$ there exists a
bound $N_a$ such that $a$ cannot be written as product of more than $N_a$
atoms.
\end{definition}

\begin{definition}
\label{def_garside_element}
For $\delta\in M$ we define the sets
$D_\delta^l = \{ x\in M : x \preceq \delta \}$
and $D_\delta^r = \{ x\in M : \delta \succeq x \}$.
The element $\delta$ is called
a {\em Garside element} of $M$ if $D_\delta^l = D_\delta^r$ and
if $D_\delta^l$ is finite and generates $M$.

The monoid $M$ is called a {\em Garside monoid} if it is atomic, has a Garside
element $\delta$ and if for all $a,b\in M$ a gcd and a lcm
of $a$ and $b$ exist.  In this case, the lcm and gcd of $a$ and $b$ are unique;
we denote them by $a\vee b$ and $a\wedge b$.  We call the elements of
$D_\delta^l$ the {\em simple elements} of $M$.
\end{definition}

\begin{theorem}
\label{thm_garside_group}
Let $M$ be a Garside monoid with Garside element $\delta$ and group
of fractions $G$.
\begin{enumerate}
\renewcommand{\labelenumi}{\textup{(\alph{enumi})}}
\item $M$ embeds into $G$.
\item If $a$ is an atom of $M$ then $a \preceq \delta$.
\item $M$ is invariant under conjugation with $\delta$.
\item For every $x\in G$ there are integers $r$ and $s$ such that
      $\delta^r \preceq x \preceq \delta^s$.
\item There is an integer $k$ such that $\delta^k$ is central in $G$.
\end{enumerate}
\end{theorem}

\begin{definition}
\label{def_garside_group}
Let $M$ be a Garside monoid. Its group of fractions $G$ is called a
{\em Garside group}.  We identify
the elements of $M$ with their images in $G$ and call them the
{\em positive elements} of $G$.  Let
$\tau\! : x\mapsto x^\delta = \delta^{-1}x\delta$ be the automorphism of $G$
induced by conjugation with $\delta$.

The partial orderings $\preceq$ and $\succeq$, and thus the notions of left gcd
and left lcm, can be extended to $G$ as follows. For $a,b\in G$, we say
$a\preceq b$ if there exists an element $c\in M$ such that $ac=b$ and
we say $a\succeq b$
if there exists an element $c\in M$ such that $a=cb$.  Clearly $\preceq$
and $\succeq$ are invariant under $\tau$.
\end{definition}

\begin{example}
\label{ex_artin_monoid}
Consider the monoid $B_n^+$ defined by the presentation
\begin{equation}
\label{eqn_artin_pres}
   \left\langle \sigma_1,\dots,\sigma_{n-1} \; \left| \;
     \begin{array}{c}
       \sigma_i \sigma_j = \sigma_j \sigma_i \quad (1\le i<j+1\le n) \\
       \sigma_i \sigma_{i+1} \sigma_i = \sigma_{i+1} \sigma_i \sigma_{i+1}
          \quad (1\le i \le n-2)
     \end{array}
   \right. \right\rangle .
\end{equation}
Its quotient group is the braid group $B_n$ on $n$ strings
\cite{artin_braid_groups}.  $B_n^+$ is a Garside monoid with Garside element
$(\sigma_1\cdots\sigma_{n-1})(\sigma_1\cdots\sigma_{n-2})\cdots(\sigma_1\sigma_2)\sigma_1$.
The positive elements of $B_n$ are simply the words in
$\sigma_1,\dots,\sigma_{n-1}$ not involving inverses of generators.
There are $n!$ simple elements, corresponding to those braids in which
any two strings cross at most once.  A simple element is described
uniquely by the permutation it induces on the strings and every permutation
of the $n$ strings corresponds to a simple element.
\end{example}

\begin{example}
\label{ex_BKL_monoid}
The monoid $BKL_n^+$ generated by $\{ a_{t,s} : n\ge t > s \ge 1 \}$
subject to the relations
\begin{equation}
\label{eqn_bkl_pres}
     \begin{array}{l}
       a_{t,s} a_{r,q} = a_{r,q} a_{t,s} 
                     \quad\mbox{if}\quad (t-r)(t-q)(s-r)(s-q) > 0 \\
       a_{t,s} a_{s,r} = a_{t,r} a_{t,s} = a_{s,r} a_{t,r}
                     \quad\mbox{if}\quad t > s > r
     \end{array}
\end{equation}
also has the braid group $B_n$ as its quotient group
\cite{braid_word_problem}.  In terms of
presentation~(\ref{eqn_artin_pres}),
$a_{t,s} = (\sigma_{t-1}\cdots\sigma_{s+1}) \sigma_s
 (\sigma_{s+1}^{-1}\cdots\sigma_{t-1}^{-1})$
is a possible choice for the generators~$a_{t,s}$.

$BKL_n^+$ is a
Garside monoid with Garside element $a_{n,n-1}a_{n-1,n-2}\cdots a_{2,1}$.
The number of simple elements of $BKL_n^+$ is $(2n)! / (n!(n+1)!)$.
Again, a simple element is described uniquely by the permutation it
induces on the strings, but not every permutation of the $n$ strings
corresponds to a simple element.
\end{example}

\begin{notation}
From now on let $M$ be a Garside monoid with Garside group~$G$,
Garside element $\delta$ and set of simple elements $D$.
\end{notation}

\subsection{Normal Forms}
\label{sect_normal_form}

\begin{definition}
\label{def_normal_form}
By Theorem \ref{thm_garside_group} there exist for every $x\in G$
 integers $r\ge0$ and $k$ such that
$\delta^k \preceq x \preceq \delta^{k+r}$.  Choose $k$ maximal and $r$
minimal subject to this condition.  We call $k$ the {\em infimum}, denoted
by $\inf(x)$, $r$ the {\em canonical length}, denoted $\len(x)$, and $k+r$
the {\em supremum}, denoted by $\sup(x)$, of $x$.

There are uniquely defined elements $A_1,\dots,A_r\in D$ such that
$x = \delta^k A_1\cdots A_r$ and $A_i^{-1}\delta\wedge A_{i+1} = 1$ for
$i=1,\dots,r-1$.  We call this representation of $x$ the {\em normal form}
of $x$.
Obviously, $A_i = \delta \wedge (\delta^k A_1\cdots A_{i-1})^{-1}x$ for
$i=1,\dots,r$.
Note that, as $A_i^{-1}\delta\preceq\delta$, we have
$A_{i+1}\cdots A_r \wedge A_i^{-1}\delta =
 A_{i+1}\cdots A_r \wedge \delta \wedge A_i^{-1}\delta =
 A_{i+1} \wedge A_i^{-1}\delta = 1$.
\end{definition}

\subsection{Super Summit Sets}
\label{sect_super_summit}
The notion of super summit sets was developed in \cite{braid_garside_alg}
and \cite{braid_positive_braids} in the context of braid groups and extended
to Garside groups in \cite{braid_conjugacy_gaussian}.  It is crucial for
testing conjugacy in Garside groups.  More details and proofs of the
results quoted in this section can be found in the references above.

\begin{definition}
\label{def_super_summit}
Let $x \in G$ and denote by $x^G$ the set of conjugates of $x$. Let
$\sinf(x) = \max\{ \inf(y) : y \in x^G \}$ and
$\ssup(x) = \min\{ \sup(y) : y \in x^G \}$.

The set
$ S_x = \{ y \in x^G : \inf(y) = \sinf(x),\, \sup(y) = \ssup(x) \} $
is called the {\em super summit set} of $x$.  We define
$\slen(x) = \ssup(x) - \sinf(x)$.
\end{definition}

\begin{definition}
\label{def_cycling_decycling}
Let $\delta^k A_1\cdots A_r \in G$ be the normal form of $x\in G$.  If
$r = 0$, let $\cyc(x) = \dec(x) = x$, otherwise let
$\cyc(x) = x^{\tau^{-k}(A_1)}$ and $\dec(x) = x^{A_r^{-1}}$.
We call $\cyc(x)$ and $\dec(x)$ the {\em cycling} of $x$ and the
{\em decycling} of $x$, respectively.
\end{definition}

\begin{theorem}
[
\cite{braid_positive_braids},
 \cite{braid_conjugacy_gaussian}, \cite{braid_infimum}]
\label{thm_super_summit}
Let $x\in G$.
\begin{enumerate}
\renewcommand{\labelenumi}{\textup{(\alph{enumi})}}
\item $S_x$ is finite and not empty.
\item A representative of $S_x$ can be obtained effectively by
 applying a finite sequence of cycling and decycling operations to $x$.
\item If $y\in S_x$ then $\cyc(y)\in S_x$ and $\dec(y)\in S_x$.
\item For all $y \in G$, $\tau(\cyc(y)) = \cyc(\tau(y))$ and
 $\tau(\dec(y)) = \dec(\tau(y))$.
\end{enumerate}
\end{theorem}

The following result is crucial for computing the super summit set of an
element.

\begin{theorem}
[El-Rifai, Morton \cite{braid_positive_braids}; Picantin \cite{braid_conjugacy_gaussian}]
\label{thm_super_summit_convex}
Let $x\in G$.
\begin{enumerate}
\renewcommand{\labelenumi}{\textup{(\alph{enumi})}}
\item For any $y,z\in S_x$ there exists $u\in M$ such that $y^u=z$.
\item If $y\in S_x$ and $u\in M$ such that $y^u\in S_x$ then
 $y^{\delta\wedge u}\in S_x$.
\item For any $y,z\in S_x$ there exist elements
 $y_0,\dots,y_t \in S_x$ and elements $c_1,\dots,c_t\in D$ such that
 $y_0 = y$, $y_t = z$ and $y_{i-1}^{c_i} = y_i$ for $i=1,\dots,t$.
\end{enumerate}
\end{theorem}

Hence $S_x$ can be computed as follows.  First obtain
$\tilde{x} \in S_x$ according to Theorem \ref{thm_super_summit} (b)
and set $S = \{ \tilde{x} \}$.
Now keep conjugating elements of $S$ by simple elements and add
those conjugates with infimum $\sinf(x)$ and supremum $\ssup(x)$ to
$S$.  When no new elements of $S_x$ can be found using this method,
that is, $S = \{ y^c : y \in S,\; c \in D,\; y^c \in S_x \}$, then $S = S_x$.

\smallskip
Franco and Gonz\'alez-Meneses improved this algorithm as follows.

\begin{theorem}[Franco, Gonz\'alez-Meneses \cite{braid_conjugacy}]
\label{thm_super_summit_gcd}
Let $x\in G$, $y\in S_x$ and $u,v\in D$. If $y^u \in S_x$ and $y^v \in S_x$
then $y^{u\wedge v} \in S_x$.
\end{theorem}


Hence, for an element $y \in S$ in the algorithm outlined above,
only the conjugates by those elements which are minimal
with respect to $\preceq$ in the set $\{ c \in D : c \ne 1,\; y^c \in S_x \}$
have to be considered.  Franco and Gonz\'alez-Meneses remark in
\cite{braid_conjugacy} that the number of such $\preceq$-minimal elements
is bounded by the number of atoms in $M$ and give an algorithm for computing
them.

\subsection{Testing Conjugacy of Elements}
\label{sect_super_summit_conjugacy}
Since $S_x$ by definition only depends on the conjugacy class of $x$, 
conjugacy of elements $x$ and $y$ of $G$ can be tested as follows
\cite{braid_positive_braids,braid_conjugacy_gaussian,braid_conjugacy}.

Compute representatives $\tilde{x}$ of $S_x$ and $\tilde{y}$ of $S_y$
according to Theorem \ref{thm_super_summit} (b).  If
$\inf(\tilde{x}) \ne \inf(\tilde{y})$ or
$\sup(\tilde{x}) \ne \sup(\tilde{y})$ then $x$ and $y$ are not conjugate.
Otherwise, start computing $S_x$ as described in Section
\ref{sect_super_summit}.  The elements $x$ and $y$ are conjugate if
and only if $\tilde{y} \in S_x$.  Note that if $x$ and $y$ are conjugate,
an element conjugating $x$ to $y$ can be found by keeping track of the
conjugations during the computations of $\tilde{x}$, $\tilde{y}$ and $S_x$.

\begin{remark}
\label{rem_super_summit_conjugacy}
It is obvious that in the worst case, both the space and the time
requirements of the algorithm outlined above are proportional to
the cardinality of $S_x$.

In the cases of the monoids $B_n^+$ and $BKL_n^+$, the only known upper
bounds for the size of $S_x$ are exponential in $n$ and $\len(x)$.
It is conjectured however, that for fixed $n$, at least for $B_n^+$
a polynomial bound in $\len(x)$ exists \cite{braid_epsteinetal}.

Nevertheless, the rapidly growing super summit sets make computations
in general infeasible for values larger than $n\approx 10$ due to lack
of memory.

Note also that distributing the computation of $S_x$ is not practical,
as the set $S$ defined in Section \ref{sect_super_summit}
is constantly accessed and modified by all nodes.
\end{remark}

\subsection{Ultra Summit Sets}
\label{sect_ultra_summit}
\begin{definition}
\label{def_ultra_summit}
By Theorem \ref{thm_super_summit}, the super summit
set $S_x$ of $x\in G$ can be made into a finite directed graph
$\Gamma_x$ with set of vertices $S_x$ and set of edges
$\{ (y,\cyc(y)) : y \in S_x \}$.  Obviously, $\tau$ induces an
automorphism of $\Gamma_x$.

Let $U_x$, the {\em ultra summit set} of $x$, be the subset of vertices
which are contained in a circuit of $\Gamma_x$, that is,
$U_x = \{ y \in S_x : \cyc^k(y) = y \mbox{ for some $k>0$} \}$.

For $y\in S_x$, define the {\em trajectory}
$T_y = \{ \cyc^k(y) : k \ge 0 \}$.  A representative of $U_x$
can be obtained by computing $T_y$ for an arbitrary $y\in S_x$.
For any $z\in T_y$, computing $s_z\in M$ satisfying $y^{s_z} = z$ is
straightforward.
\end{definition}

\smallskip
The following main result of this paper will be proved in
Section \ref{sect_proof_main}.

\begin{theorem}
\label{thm_ultra_summit_gcd}
Let $x\in G$, $y \in U_x$ and let $u,v\in M$ such that $y^u \in U_x$ and
$y^v \in U_x$.  Then $y^{u\wedge v} \in U_x$.
\end{theorem}

\begin{corollary}
\label{coro_ultra_summit_convex}
Let $x\in G$ and $y,z\in U_x$. There exist elements
$y_0,\dots,y_t \in U_x$ and elements $c_1,\dots,c_t\in D$ such that
$y_0 = y$, $y_t = z$ and $y_{i-1}^{c_i} = y_i$ for $i=1,\dots,t$.
\end{corollary}
\proof
We may assume $y\ne z$. First note that $y\in U_x$ implies
$y^\delta = \tau(y) \in U_x$ as
$\tau$ is an automorphism of $\Gamma_x$.
By Theorem~\ref{thm_super_summit_convex} (a), there exists $u \in M$
with $y^u=z$.  Let $s = \sup(u)$.  Choose
$c_1 = \delta \wedge u \in D$ and let $y_1 = y^{c_1}$ and
$\tilde{u} = c_1^{-1}u \in M$.  By
Theorem~\ref{thm_ultra_summit_gcd}, $y_1 \in U_x$.  Moreover,
$y_1^{\tilde{u}} = z$ and $\sup(\tilde{u}) < s$.  Iteration yields
$y_1,\dots,y_t \in U_x$ and $c_1,\dots,c_t\in D$ as desired.
\qed

\begin{definition}
\label{def_minimal_simple_elements}
Let $x\in G$ and $y\in U_x$.
\begin{enumerate}
\renewcommand{\labelenumi}{\textup{(\alph{enumi})}}
\item 
 For any $s\in D$, Theorem
 \ref{thm_ultra_summit_gcd} implies the existence
 of a unique $\preceq$-minimal element $c_s = c_s(y)$ satisfying
 $s\preceq c_s \preceq \delta$ and $y^{c_s} \in U_x$.

\item
 Define $D_y = \{ u \in D \setminus\{1\} : y^u \in U_x \}$ and let
 $C_y$ be the set of elements of $D_y$ which are $\preceq$-minimal
 in $D_y$. Clearly $C_y \subseteq \{ c_a(y) : a \in A \}$, where
 $A$ is the set of atoms of $M$. In particular, $|C_y| \le |A|$.
\end{enumerate}
\end{definition}

\begin{corollary}
\label{coro_minimal_simple_elements}
Let $x\in G$ and $U\subseteq U_x$.
If $\{ y^c : y \in U,\; c \in C_y \}\subseteq U$ then $U = U_x$.
\end{corollary}
\proof
This follows directly from Corollary \ref{coro_ultra_summit_convex}.
\qed

\smallskip
The following result will also be proved in Section \ref{sect_proof_main}.

\begin{theorem}
\label{thm_minimal_simple_elements}
Let $x\in G$, $y\in U_x$ and $z\in T_y$.  For any $s\in C_z$ there exists
$t\in C_y$ such that $z^s\in T_{y^t}$.
\end{theorem}

\begin{algo}
\label{alg_compute_ultra_summit}
Given an element $x$ of a Garside group, the following algorithm
computes the ultra summit set $U_x$ of $x$.
\upshape
\begin{algorithmic}
\STATE  Compute $\tilde{x} \in U_x$, set $U = T_{\tilde{x}}$ and
        $U_0 = \emptyset$.
\IF{$\tilde{x} = \delta^k$ for some $k$} 
   \STATE \textbf{return} $\{ \delta^k \}$
\ENDIF
\WHILE{$U \ne U_0$}
   \STATE Let $y_1,\dots,y_m \in U$ such that
          $U = U_0 \cup T_{y_1} \cup \dots \cup T_{y_m}$. Set $U_0 = U$.
   \FOR{$y \in \{y_1,\dots,y_m\}$}
      \STATE Compute $C_y$ and set $U = U \,\cup\, \bigcup_{c \in C_y}
             T_{y^c}$. \hfill$[\ast]$
   \ENDFOR
\ENDWHILE
\STATE \textbf{return} $U$
\end{algorithmic}
The computation of the set $C_y$ in step $[\ast]$ will be discussed
in Section \ref{sect_compute_minimal}.
\end{algo}

Two elements $x$ and $y$ of $G$ are conjugate in $G$ if and only if
$U_x = U_y$, or indeed, if and only if $U_x \cap U_y \ne \emptyset$.
Hence, conjugacy of elements $x$ and $y$ of $G$ can
be tested, and a conjugating element can be computed, as outlined in
Section \ref{sect_super_summit_conjugacy}, using ultra summit sets
instead of super summit sets.

\newcommand{\refthmultrasummitgcd}{\ref{thm_ultra_summit_gcd}}
\newcommand{\refthmminimalsimpleelements}{\ref{thm_minimal_simple_elements}}
\section{Proof of Theorems \refthmultrasummitgcd{} and
 \refthmminimalsimpleelements}
\label{sect_proof_main}
Throughout this section let $x\in G$ be an element of its super summit
set with non-zero canonical length, that is, let $\delta^k A_1\cdots A_r$
be the normal form of $x$, with $r>0$,
$k = \inf(x) = \sinf(x)$ and $r+k = \sup(x) = \ssup(x)$.

We need to understand how the normal forms of conjugates of $x$ are related
to the normal form of $x$.

\begin{proposition}
\label{prop_moving_constituents}
Let $x$ be as above and let $u\in M$ such that $x^u\in S_x$.
There are elements $u_0,\dots,u_r$ in $M$ such that $u_0 = \tau^k(u)$,
$u_r = u$ and the normal form of $x^u$ is
$\delta^k (u_0^{-1}A_1u_1)\cdots(u_{r-1}^{-1}A_r u_r)$.
Here, the factors in brackets are understood to be the simple elements
occurring in the normal form of $x^u$. Explicitly,
$u_i = A_{i+1}\cdots A_r u \,\wedge\, \delta \tau(A_i^{-1}u_{i-1})$.
\end{proposition}
\proof
Let $u_0 = \tau^k(u)$ and $u_r = u$. Define
$w_1 = \delta^{-k}x^u = u_0^{-1}A_1\cdots A_r u_r$ and
 $w_{i+1} = (w_i\wedge\delta)^{-1} w_i$ for $i=1,\dots,r-1$. 
By the observation in Definition \ref{def_normal_form}, $w_i$ has infimum $0$
and canonical length $r+1-i$ and the normal
form of $x^u$ is $\delta^k (w_1\wedge\delta)\cdots(w_r\wedge\delta)$.
Assume $u_{i-1}\in M$ has been found such that
$w_i = u_{i-1}^{-1}A_i\cdots A_r u_r$.
Then, $A_i \preceq \delta \preceq u_{i-1}\delta$ implies
$u_{i-1}^{-1}A_i \preceq w_i \wedge \delta$,
that is, there is an element $u_i\in M$ such that
$w_i \wedge \delta = u_{i-1}^{-1}A_i u_i$. Now
$w_{i+1} = (w_i \wedge \delta)^{-1} w_i = u_i^{-1}A_{i+1}\cdots A_r u_r$
and $u_i = A_i^{-1}u_{i-1}(w_i\wedge\delta) =
 A_{i+1}\cdots A_r u \,\wedge\, \delta \tau(A_i^{-1}u_{i-1})$ as claimed.
\qed

\begin{corollary}
\label{coro_moving_constituents}
Let $x$ be as above and let $u, v\in M$ such that $x^u\in S_x$ and
$x^v\in S_x$.  Let $u_0,\dots,u_r$ and $v_0,\dots,v_r$ be the
positive elements obtained by applying Proposition
\ref{prop_moving_constituents} to $(x,u)$ and $(x,v)$, respectively.
\begin{enumerate}
\renewcommand{\labelenumi}{\textup{(\alph{enumi})}}
\item If $u = \delta$ then $u_i = \delta $ for $i=0,\dots,r$.
\item If $u\preceq v$ then $u_i\preceq v_i$ for $i=0,\dots,r$.  More
 specifically, if $v = u w$ with $w\in M$ and $w_0,\dots,w_r$ are the
 positive elements obtained by applying Proposition
 \ref{prop_moving_constituents} to $(x^u,w)$ then $v_i = u_i w_i$
 for $i=0,\dots,r$.
\item If $\sup(u) = b$ then $\sup(u_i) \le b$ for $i=0,\dots,r$.  In
 particular, if $u$ is simple then $u_i$ is simple for $i=0,\dots,r$.
\item If $u\wedge v = 1$ then $u_i\wedge v_i = 1$ for $i=0,\dots,r$.
\item Let $t = u\,\wedge\, v$ and let $t_0,\dots,t_r$ be the
 positive elements obtained by applying Proposition
 \ref{prop_moving_constituents} to $(x,t)$. Then $t_i = u_i \wedge v_i$
 for $i=0,\dots,r$.
\end{enumerate}
\end{corollary}
\proof
(a) By Proposition \ref{prop_moving_constituents},
$u_i = \delta \tau(A_{i+1}\cdots A_r \,\wedge\, A_i^{-1}u_{i-1})$.
As $u_0 = \delta$ and
$A_{i+1}\cdots A_r \,\wedge\, A_i^{-1}\delta = 1$ by Definition
\ref{def_normal_form}, $u_i = \delta$ follows by induction.

(b) $v_0 = u_0 w_0$ is obvious. Assume
$v_{i-1} = u_{i-1} w_{i-1}$. By Proposition \ref{prop_moving_constituents},
$w_i = (u_i^{-1}A_{i+1}u_{i+1})\cdots(u_{r-1}^{-1}A_r u_r)w \,\wedge\,
 \delta \tau\left((u_{i-1}^{-1}A_i u_i)^{-1}w_{i-1}\right)$,
whence $u_i w_i = A_{i+1}\cdots A_r v\,\wedge\,\delta\tau(A_i^{-1}v_{i-1})
 = v_i$, again using Proposition \ref{prop_moving_constituents}. Hence
the claim follows by induction.

(c) Follows from parts (a) and (b), as $\sup(u) \le b$
if and only if $u \preceq \delta^b$.

(d) $u_0\wedge v_0 = 1$ is obvious. Assume
$u_{i-1}\wedge v_{i-1} = 1$. By Proposition \ref{prop_moving_constituents},
$u_i \wedge v_i = A_{i+1}\cdots A_r (u\wedge v) \wedge A_i^{-1}\delta 
\tau(u_{i-1}\wedge v_{i-1}) =
 A_{i+1}\cdots A_r \,\wedge\, A_i^{-1}\delta = 1$, where in the last step
Definition \ref{def_normal_form} was used.  Hence the claim follows
by induction. 

(e) Note that $x^t \in S_x$ by Theorems
\ref{thm_super_summit_convex} (b) and \ref{thm_super_summit_gcd}, that is,
Proposition \ref{prop_moving_constituents} can be applied to $(x,t)$.
The claim then follows from parts (b) and (d), writing $u=t\bar{u}$ and
$v=t\bar{v}$ with $\bar{u}\wedge\bar{v} = 1$.
\qed

\begin{lemma}
\label{lem_transport}
Let $x$ be as above, $u\in M$ such that $x^u \in S_x$.  Let $u_0,\dots,u_r$
be the positive elements obtained by applying Proposition
\ref{prop_moving_constituents} to $(x,u)$.  Let $\phi_x(u) = \tau^{-k}(u_1)$.
\begin{enumerate}
\renewcommand{\labelenumi}{\textup{(\alph{enumi})}}
\item $\phi_x(u) \in M$ satisfies
 $\cyc(x^u) = \cyc(x)^{\phi_x(u)}$.
\item $\sup(\phi_x(u)) \le \sup(u)$. In particular, if $u$ is simple then
 $\phi_x(u)$ is simple.
\item The conjugating element along any path in the diagram
\[
\xymatrix{
x^u \ar[rr]|-{\;\cyc\;}^(.45){\rule[-1ex]{0ex}{1ex}\tau^{-k}(u_0^{-1}A_1u_1)}
  & &
  \cyc(x^u) \\
x \ar[u]^{u} \ar[rr]|-{\;\cyc\;}^(.45){\rule[-1ex]{0ex}{1ex}\tau^{-k}(A_1)} & &
  \cyc(x) \ar[u]_{\phi_x(u)}
}
\]
only depends on the starting point and the end point of the path.
\end{enumerate}
\end{lemma}
\proof
(a) follows from $\cyc(x)^{\tau^{-k}(u_1)} = x^{\tau^{-k}(A_1 u_1)} =
 (x^u)^{\tau^{-k}(u_0^{-1}A_1 u_1)} = \cyc(x^u)$.
The conjugating element along the circuit $x \rightarrow
 x^u \rightarrow \cyc(x^u) \rightarrow \cyc(x) \rightarrow x$ is
 $u\cdot \tau^{-k}(u_0^{-1}A_1u_1) \cdot \phi_x(u)^{-1} \cdot
 \tau^{-k}(A_1)^{-1} = 1$, proving (c).
Part (b) follows from Corollary \ref{coro_moving_constituents} (c).
\qed

\begin{definition}
\label{def_transport}
In the situation of Lemma \ref{lem_transport}, we call
$\phi_x(u)$ the {\em transport} of $u$ along $x \rightarrow \cyc(x)$.
If $x$ is obvious from the context, we define $u^{(0)} = u$ and
$u^{(i+1)} = \phi_{\textup{\bf\scriptsize c}^{i}(x)}(u^{(i)})$ for $i\ge 0$.
\end{definition}

\begin{lemma}
\label{lem_transport_injective}
Let $x$ be as above and let $u, v\in M$ such that $x^u = x^v \in S_x$.
If $\phi_x(u) = \phi_x(v)$ then $u = v$.
\end{lemma}
\proof
Let $u_0,\dots,u_r$ and $v_0,\dots,v_r$ be the
positive elements obtained by applying Proposition
\ref{prop_moving_constituents} to $(x,u)$ and $(x,v)$, respectively.
As $x^u = x^v$, we have $(u_0^{-1}A_1u_1)
  = \delta \wedge \delta^{-k}x^u = \delta \wedge \delta^{-k}x^v
  = (v_0^{-1}A_1v_1)$.  The claim then follows from
Lemma \ref{lem_transport} (c).
\qed

\begin{lemma}
\label{lem_transport_periodic}
Let $x$ be as above, let $u\in M$ such that $x^u \in S_x$ and let
$\cyc^N(x) = x$ and $\cyc^N(x^u) = x^u$ for some integer $N>0$.
There is an integer $m>0$ such that $u^{(mN)} = u$, where we use the notation
from Definition~\ref{def_transport}.
\end{lemma}
\proof
By Lemma \ref{lem_transport} (b), $u^{(iN)}\in M$ and
$\sup(u^{(iN)})\le\sup(u)$ for every integer $i\ge 0$.  Since the number
of such elements is at most $|D|^{\sup(u)}$, in particular finite, there
must exist integers $i_2 > i_1 \ge 0$ such that $u^{(i_1N)} = u^{(i_2N)}$;
let $i_2$ be minimal subject to this condition.  Assume $i_1>0$. Then we
can for $l=1,\dots,N$ conclude $u^{(i_1N-l)} = u^{(i_2N-l)}$ from
\begin{eqnarray*}
  && \phi_{\textup{\bf\scriptsize c}^{(N-l)}(x)}\left(u^{(i_1N-l)}\right) =
     \phi_{\textup{\bf\scriptsize c}^{(i_1N-l)}(x)}\left(u^{(i_1N-l)}\right) = u^{(i_1N-(l-1))} \\
  && = u^{(i_2 N-(l-1))} = \phi_{\textup{\bf\scriptsize c}^{(i_2N-l)}(x)}\left(u^{(i_2N-l)}\right)
     = \phi_{\textup{\bf\scriptsize c}^{(N-l)}(x)}\left(u^{(i_2N-l)}\right)  \;,
\end{eqnarray*}
using Lemma \ref{lem_transport_injective}. In particular,
$u^{((i_1-1)N)} = u^{((i_2-1)N)}$, contradicting the minimality of $i_2$.
Hence, $i_1=0$ and $u^{(i_2N)} = u^{(0)} = u$ as claimed.
\qed

\begin{theorem}
\label{thm_ultra_summit_trivial_gcd}
Let $x$ be as above, $u, v\in M$ such that $u\wedge v = 1$.
If $x^u \in U_x$ and $x^v \in U_x$ then $x\in U_x$.
\end{theorem}
\proof
First note that we may assume that $\cyc(x)\in U_x$, since if $x$ is a
counterexample with $\cyc(x)\notin U_x$, consider $\bar{x} = \cyc(x)\in S_x$,
$\bar{u} = \phi_x(u)$ and $\bar{v} = \phi_x(v)$.  Clearly,
$\bar{x}^{\bar{u}} = \cyc(x^u)\in U_x$ and
$\bar{x}^{\bar{v}} = \cyc(x^v)\in U_x$.
Moreover, $\bar{u}\wedge\bar{v}=1$ by Corollary
\ref{coro_moving_constituents} (d).  Repeating this process finitely many
times, we arrive at a counterexample $x$ with $\cyc(x)\in U_x$.

Choose $N>0$ such that $\cyc^N(x^u) = x^u$, $\cyc^N(x^v) = x^v$, and
$\cyc^{N+1}(x) = \cyc(x)$.  We use the notation from Definition
\ref{def_transport}. According to Lemma \ref{lem_transport_periodic}, we
can further assume that $u^{(N+1)} = u^{(1)}$ and $v^{(N+1)} = v^{(1)}$,
replacing $N$ by a suitable multiple if necessary.
Now consider the conjugations by the conjugating elements indicated in the
following diagram.
\[
\xymatrix{
x^u \ar[r]|-{\;\cyc\;}^(.45){\rule[-1ex]{0ex}{1ex}\alpha_u} &
  \cyc(x^u) \ar[r]|-{\;\cyc\;} & \cdots \ar[r]|-{\;\cyc\;} &
  **{!<-1.45em,0ex>}\cyc^N(x^u) = x^u
            \ar[r]|-{\;\cyc\;}^(.65){\rule[-1ex]{0ex}{1ex}\beta_u} &
  \cyc(x^u)
  \\
x \ar[u]^{u} \ar[d]_{v} \ar[r]|-{\;\cyc\;}^(.45){\rule[-1ex]{0ex}{1ex}\alpha} &
  \cyc(x) \ar[u]_{u^{(1)}} \ar[d]^{v^{(1)}} \ar[r]|-{\;\cyc\;}
  & \cdots \ar[r]|-{\;\cyc\;} &
  \cyc^N(x) \ar[u]_{u^{(N)}} \ar[d]^{v^{(N)}}
            \ar[r]|-{\;\cyc\;}^(.5){\rule[-1ex]{0ex}{1ex}\beta} &
  **{!<-2em,0ex>}\cyc^{N+1}(x) = \cyc(x) \ar[u]_{u^{(N+1)}} \ar[d]^{v^{(N+1)}}
  \\
x^v \ar[r]|-{\;\cyc\;}^(.45){\rule[-1ex]{0ex}{1ex}\alpha_v} &
  \cyc(x^v) \ar[r]|-{\;\cyc\;} & \cdots \ar[r]|-{\;\cyc\;} &
  **{!<-1.45em,0ex>}\cyc^N(x^v) = x^v
              \ar[r]|-{\;\cyc\;}^(.65){\rule[-1ex]{0ex}{1ex}\beta_v} &
  \cyc(x^v)
}
\]
Obviously, $\alpha_u = \tau^{-k}(\delta \wedge \delta^{-k}x^u) = \beta_u$ and
$\alpha_v = \tau^{-k}(\delta \wedge \delta^{-k}x^v) = \beta_v$ and by Corollary
\ref{coro_moving_constituents} (d), we have $u^{(i)}\wedge v^{(i)}$
for $i=1,\dots,N$.  
Hence,
\begin{eqnarray*}
  \alpha^{-1} & = & \alpha^{-1} (u \,\wedge\, v)
    = \alpha^{-1} u \,\wedge\, \alpha^{-1} v
    =  u^{(1)}\alpha_u^{-1} \,\wedge\, v^{(1)}\alpha_v^{-1} \\
              & = & u^{(N+1)}\beta_u^{-1} \,\wedge\, v^{(N+1)}\beta_v^{-1}
    = \beta^{-1} u^{(N)} \,\wedge\, \beta^{-1} v^{(N)} = \beta^{-1}
\end{eqnarray*}
where we used Lemma \ref{lem_transport} (c) four times.
We conclude $x\in U_x$ from
\[
 x = \cyc(x)^{\alpha^{-1}} = \cyc(x)^{\beta^{-1}}
 = \left(\cyc^{N+1}(x)\right)^{\beta^{-1}} = \cyc^N(x) \; .\qed
\]

\pagebreak
Theorems \refthmultrasummitgcd{} and \refthmminimalsimpleelements{} now
follow easily.

\begin{oldtheorem}{\refthmultrasummitgcd}
Let $x\in G$, $y \in U_x$ and let $u,v\in M$ such that $y^u \in U_x$ and
$y^v \in U_x$.  Then $y^{u\wedge v} \in U_x$.
\end{oldtheorem}
\proof
If $\sinf(x) = \ssup(x) = k$ then $U_x = S_x = \{ \delta^k \}$ and
the claim follows from Theorems \ref{thm_super_summit_convex} (b) and
\ref{thm_super_summit_gcd}. Hence assume $\ssup(x) > \sinf(x)$.

Let $t = u\wedge v$. Then $u=t\bar{u}$, $v=t\bar{v}$ with
$\bar{u}\wedge \bar{v} = 1$.  By Theorems
\ref{thm_super_summit_convex} (b) and \ref{thm_super_summit_gcd}, $y^t\in S_x$.
As $(y^t)^{\bar{u}} = y^u\in U_x$ and $(y^t)^{\bar{v}} = y^v\in U_x$,
Theorem \ref{thm_ultra_summit_trivial_gcd} implies $y^t\in U_x$.
\qed

\begin{oldtheorem}{\refthmminimalsimpleelements}
Let $x\in G$, $y\in U_x$ and $z\in T_y$.  For any $s\in C_z$ there exists
$t\in C_y$ such that $z^s\in T_{y^t}$.
\end{oldtheorem}
\proof
Consider the restriction
$\phi = \phi_y{}_{|_{D_y \cup \{1\}}} : D_y \cup \{1\}
 \rightarrow D_{\textup{\bf\scriptsize c}(y)} \cup \{1\}$
of $\phi_y$ to $D_y \cup \{1\}$. By Lemma \ref{lem_transport_periodic},
$\phi$ is bijective and, using Corollary \ref{coro_moving_constituents} (b),
$\phi(u)\preceq\phi(v)$ if and only if $u\preceq v$ holds for all
$u,v\in D_y \cup \{1\}$.  We hence obtain
$C_{\textup{\bf\scriptsize c}(y)} = \{ \phi_x(u) : u \in C_y \}$
from which the claim follows by induction.
\qed

\section{A Probabilistic Approach to the Conjugacy Search Problem}
\label{sect_probabilistic}

Given elements $x,y\in G$ which are conjugate in $G$, we can use the
structure of the graph $\Gamma_x$ for computing an element $s\in G$
satisfying $x^s=y$ without having to compute the entire ultra summit set $U_x$.

Applying cycling and decycling operations to $x$ and $y$, respectively,
we can obtain $\tilde{x},\tilde{y}\in U_x = U_y$ as well as $s_x, s_y\in G$
satisfying $x^{s_x} = \tilde{x}$ and
$y^{s_y} = \tilde{y}$. For $z\in T_{\tilde{x}}$, that is,
$z = \cyc^k(\tilde{x})$ for some $k$, let $s(z)$ satisfy
$\tilde{x}^{s(z)} = z$.

\begin{algo}
\label{alg_conjugacy_search}
Given a Garside group $G$ and elements $x, y\in G$ which are
conjugate in $G$, the following Las Vegas algorithm computes
an element $s\in G$ such that $x^s=y$.
\upshape
\begin{algorithmic}
\STATE  Compute $\tilde{x}$, $s_x$, $T_{\tilde{x}}$ and
           $\{ s(z) : z \in T_{\tilde{x}} \}$ as above.
   \STATE  Compute $\tilde{y}$ and $s_y$ as above.\ Set
              $z = \tilde{y}$ and $s = s_y$.
\LOOP
   \IF{$z\in T_{\tilde{x}}$}
      \STATE \textbf{return} $s_x\cdot s(z)\cdot s^{-1}$
   \ENDIF
   \STATE  Choose a random atom $a$ of $M$.\ Compute $c_a = c_a(z)$.
           \hfill$[\ast]$
   \STATE  Set $z = z^{c_a}$, $s = s\cdot c_a$.
\ENDLOOP
\end{algorithmic}
The computation of $c_a$ in step $[\ast]$ (recall Definition
\ref{def_minimal_simple_elements}) will be discussed in
Section \ref{sect_compute_minimal}.
\end{algo}

\begin{remark}
\label{rem_ultra_summit_conjugacy}
The expected number of passes through the loop in Algorithm
\ref{alg_conjugacy_search} is the number of circuits of the graph
$\Gamma_x$.  This loop can easily be parallelised, since no
communication between nodes is necessary.
\end{remark}

\section{Computing Minimal Elements}
\label{sect_compute_minimal}
Throughout this section let $x\in G$ be an element of its ultra summit
set with normal form $\delta^k A_1\cdots A_r$, where $r>0$, and let
$N$ be the minimal positive integer satisfying $\cyc^N(x) = x$.

In this section we show how the elements $c_s = c_s(x)$ ($s\in D$) and
the set $C_x$ introduced in Definition~\ref{def_minimal_simple_elements}
can be computed efficiently.

For any $s\in D$, Theorem \ref{thm_super_summit_gcd} implies the existence
of a unique $\preceq$-minimal element $\rho_s = \rho_s(x)$ satisfying
$s\preceq \rho_s\preceq\delta$ and $x^{\rho_s}\in S_x$.  An algorithm
for computing $\rho_s$ is given in \cite{braid_conjugacy}.

Note that $\rho_s\preceq c_s$ since $U_x\subseteq S_x$. Moreover,
if $s = 1$ then $c_s = \rho_s = 1$.

\begin{definition}
\label{def_fixed_points}
Let $u\in D$ such that $x^u\in S_x$. Using the notation from
Definition \ref{def_transport}, we consider the elements $u^{(iN)}$
($i=0,1,\dots$).  By Lemma \ref{lem_transport} (b) and since $D$ is finite,
there are integers $i_2 > i_1 \ge 0$ such that $u^{(i_1N)} = u^{(i_2N)}$.
Let $i_1$ and $i_2$ be minimal subject to this condition and define $l_x(u) = i_2 - i_1$
and $F_x(u) = \{ u^{(iN)} : i_1\le i < i_2 \}$.

Note that $1 \in F_x(u)$ if and only if $F_x(u) = \{ 1 \}$.
Moreover, if $x^u\in U_x$ then $i_1=0$ by Lemma
\ref{lem_transport_periodic}, that is, $u\in F_x(u)$.

\end{definition}

\begin{lemma}
\label{lem_fixed_points}
Let $u\in D$ such that $x^u \in S_x$, let $v \in F_x(u)$ and let $l=l_x(u)$.
Then, $v^{(ilN)} = v$ for
all integers $i>0$. Moreover, $x^v \in U_x$.
\end{lemma}
\proof
As $v^{(lN)} = v$, the first claim follows by induction. For the second
claim note that
$\cyc^{lN}(x^v) = x^{(v^{(lN)})} = x^v$, whence $x^v \in U_x$.
\qed

\begin{lemma}
\label{lem_no_proper_prefix}
Let $s \in D$. If $c_s \preceq c_s^{(iN)}$ for some $i>0$
then $c_s^{(iN)} = c_s$.
\end{lemma}
\proof
Let $c_s^{(iN)} = c_s \gamma$ with $\gamma\in M$. By induction,
$c_s \gamma \preceq c_s^{(\beta iN)}$ for all $\beta \ge 1$ from
Corollary \ref{coro_moving_constituents} (b).  Using
Lemma \ref{lem_fixed_points}, this in particular implies
$c_s \preceq c_s \gamma \preceq c_s^{(l_x(c_s) iN)} = c_s$,
that is, $\gamma = 1$.
\qed

\begin{lemma}
\label{lem_fixed_points_of_rho}
Let $p,s \in D$ satisfy $p\preceq c_s$ and $x^p \in S_x$.  Let $F = F_x(p)$.
\begin{enumerate}
\renewcommand{\labelenumi}{\textup{(\alph{enumi})}}
\item If there exists $v \in F$ such that $s\preceq v$ then $c_s = v$.
\item If $F \ne \{ 1 \}$ and $s\not\preceq v$ for all $v\in F$ then
      $c_s$ is not $\preceq$-minimal in $D_x$.
\end{enumerate}
\end{lemma}
\proof
First note that by Corollary \ref{coro_moving_constituents} (b),
$p^{(i)} \preceq c_s^{(i)}$ for all $i>0$.

(a) As $s\preceq v$ and $x^v\in U_x$ by Lemma \ref{lem_fixed_points},
 minimality of $c_s$ implies $c_s\preceq v$. Now $v = p^{(iN)}$ for some $i$,
 whence $c_s\preceq v = p^{(iN)} \preceq c_s^{(iN)}$.  Lemma
 \ref{lem_no_proper_prefix} yields $v = c_s$.

(b) Let $i$ be a multiple of $l_x(c_s)$ sufficiently large so that
 $v = p^{(iN)} \in F$.  Since $1\notin F$, we have $v\in D_x$ by Lemma
 \ref{lem_fixed_points} and Corollary \ref{coro_moving_constituents} (c).
 Moreover, again using Lemma \ref{lem_fixed_points},
 $v = p^{(iN)} \preceq c_s^{(iN)} = c_s$ and $v\ne c_s$, since
 $s\not\preceq v$.
\qed

\begin{lemma}
\label{lem_trivial_implies_prefix}
Let $p,s \in D \setminus \{ 1 \}$ such that $x^p \in S_x$. If there exists
an integer $i>0$ such
that $p^{(i)} = 1$ then $p \wedge \tau^{-k}(A_1) \ne 1$.

If moreover $p\preceq c_s$ and $c_s \not\preceq \tau^{-k}(A_1)$ then
$c_s$ is not $\preceq$-minimal in $D_x$.
\end{lemma}
\proof
If $p^{(1)}=1$ then Proposition \ref{prop_moving_constituents} implies
$\tau^k(p) \preceq A_1$. Thus we assume $p^{(1)}\ne 1$ and $i>1$.
Let $\delta^k B_1\cdots B_r$ be the normal form of
$\cyc(x) = x^{\tau^{-k}(A_1)}$. According to Proposition
\ref{prop_moving_constituents},
$(\tau^{-k}(A_1))^{(1)} = \phi_x(\tau^{-k}(A_1)) = \tau^{-k}(B_1)$.
By induction $(p^{(1)})^{(i-1)} = p^{(i)} = 1$ yields
\[ \left(p \wedge \tau^{-k}(A_1)\right)^{(1)}
    = p^{(1)} \wedge \left(\tau^{-k}(A_1)\right)^{(1)}
    = p^{(1)} \wedge \tau^{-k}(B_1) \ne 1 \]
using Corollary \ref{coro_moving_constituents} (e). This completes the proof
of the first claim.

Let $c = c_s \wedge \tau^{-k}(A_1) \preceq c_s$. If
$c_s \not\preceq \tau^{-k}(A_1)$ then $c \ne c_s$.
Now $p\preceq c_s$ implies $c \ne 1$ and $c\in D_x$ by
Theorem \ref{thm_ultra_summit_gcd}, since
$\cyc(x) = x^{\tau^{-k}(A_1)} \in U_x$.
\qed

\begin{definition}
\label{def_pullback}
Let $s\in D$ and let $y\in U_x$. By Theorems
\ref{thm_super_summit_convex} (b) and 
\ref{thm_super_summit_gcd} and Corollary \ref{coro_moving_constituents} (a)
and (e), there exists a unique $\preceq$-minimal element $\pi_y(s) \in D$
satisfying $y^{\pi_y(s)} \in S_x$ and $s \preceq \phi_y(\pi_y(s))$.  We call
$\pi_y(s)$ the {\em pullback} of $s$ along $y \rightarrow \cyc(y)$.
If $y$ is obvious from the context, we define $s_{(0)} = s$ and
$s_{(i+1)} = \pi_{\textup{\bf\scriptsize c}^\alpha(y)}(s_{(i)})$
for $i\ge 0$, where $0 \le \alpha\equiv -i\; (\mbox{mod}\; N)$.
\end{definition}

\begin{proposition}
\label{prop_pullback}
Let $s\in D$ and let $\delta^k B_1\cdots B_r$ be the normal form of
$y\in U_x$.  Define $b_0 = 1 \vee \tau^{-k}(B_1)s\delta^{-1}$,
$b_1 = \tau^k(s)$, $b_i = 1 \vee B_i^{-1}b_{i-1}$ for $i=2,\dots,r$
and $b = b_0 \vee b_r$.  Then $b\in D$ and $\rho_b = \pi_y(s)$.
\end{proposition}
\proof
Firstly, $\tau(b_0) = \tau^{-k}(B_1^{-1}\delta)^{-1} \cdot
(\tau^{-k}(B_1^{-1}\delta) \vee s) \in D$. Moreover, $b_1 \in D$ and
$b_i = B_i^{-1} \cdot (B_i \vee b_{i-1}) \in D$ for $i=2,\dots,r$
by induction. Hence $b\in D$.

Let $u\in D$ such that $y^u \in S_x$ and $s\preceq \phi_y(u)$.
Using the notation from Proposition \ref{prop_moving_constituents} we have
$\tau^k(s) \preceq u_1 \preceq B_1^{-1}u_0\delta =
 B_1^{-1}\delta\tau^{k+1}(u)$, that is,
$\tau^{-k}(B_1)s\delta^{-1} \preceq u$, whence $b_0 \preceq u$, as $u \in M$.
Moreover, $u_{i-1} \preceq B_i u_i$, that is, $B_i^{-1}u_{i-1} \preceq u_i$
for $i=2,\dots,r$. As $b_1 = \tau^k(s) \preceq u_1$ and $u_i \in M$
for $i=2,\dots,r$, we obtain
$b_i \preceq u_i$ for $i=1,\dots,r$ by induction; in particular,
$b_r \preceq u_r = u$.  Hence $b \preceq u$ and $\rho_b \preceq u$
by minimality of $\rho_b$, proving $\rho_b \preceq \pi_y(s)$.

Conversely, let $\rho = \rho_b$ and use again the notation from
Proposition \ref{prop_moving_constituents}.
Since $\tau^{-k}(B_1)s\delta^{-1} \preceq b_0 \preceq \rho$, we have
$\tau^k(s) \preceq B_1^{-1}\tau^k(\rho)\delta$.  On the other hand,
$\tau^k(s) = b_1 \preceq B_2b_2 \preceq \dots \preceq B_2\cdots B_rb_r
\preceq B_2\cdots B_r\rho$. Together these imply $\tau^k(s) \preceq
 B_1^{-1}\tau^k(\rho)\delta \,\wedge\, B_2\cdots B_r\rho = \rho_1$,
that is, $s \preceq \phi_y(\rho)$.
As $y^\rho \in S_x$, this proves $\pi_y(s) \preceq \rho$.
\qed

\begin{proposition}
\label{prop_pullback_yields_stable_transport}
Let $s\in D$ and consider for $i=0,1,\dots$ the elements
$s_{(iN)}$ obtained by applying Definition \ref{def_pullback} for
$y=\cyc^{N-1}(x)$.  As $D$ is finite,
there are integers $i_2 > i_1 \ge 0$ such that $s_{(i_1N)} = s_{(i_2N)}$.
Choose minimal values for $i_1$ and $i_2$, let $l = i_2-i_1$ and choose
an integer $j$ such that $jl \ge i_1$. Finally, let
$p = p_x(s) = s_{(jlN)}$.

Then, $p \preceq c_s$ and there exists $v\in F_x(p)$ with $s \preceq v$.
\end{proposition}
\proof
Let $\beta \ge j$ be a multiple of $l_x(c_s)$ large enough such
that $p^{(\beta lN)} \in F_x(p)$.
By Definition \ref{def_pullback} and Corollary
\ref{coro_moving_constituents} (b), $p = s_{(jlN)} = s_{(\beta lN)}$
is the unique \linebreak
$\preceq$-minimal
element satisfying $x^p \in S_x$ and $s \preceq p^{(\beta lN)}$.  Since
$x^{c_s} \in U_x\subseteq S_x$ and $s \preceq c_s = c_s^{(\beta lN)}$,
we obtain $p \preceq c_s$.
\qed

\pagebreak

\begin{algo}
\label{alg_compute_minimal}
Given $s\in D$ and a boolean value \texttt{\textup{f}} indicating
whether elements which are known not to be $\preceq$-minimal in $D_x$
should be discarded, the following algorithm returns $c_s$ or identifies
it as not $\preceq$-minimal in $D_x$.
\begin{algorithmic}
\upshape
\STATE  Compute $\rho_s$ as described in \cite{braid_conjugacy} and
        compute $F_x(\rho_s)$.
\IF{$\exists\, v \in F_x(\rho_s)$ such that $s \preceq v$}
   \STATE \textbf{return} $v$
\ENDIF
\IF{\texttt{f} {\bf{}and} $F_x(\rho_s) \ne \{1\}$}
\STATE \textbf{return} \texttt{not minimal}
\ENDIF
\STATE  Compute $p_x(s)$ and $F_x(p_x(s))$. \hfill$[\ast]$
\STATE  Choose $v \in F_x(p_x(s))$ such that $s \preceq v$.
\STATE  \textbf{return} $v$
\end{algorithmic}
In the case that \texttt{\textup{f}} is \texttt{\textup{true}}, the algorithm
can be aborted returning \texttt{\textup{not minimal}} in step $[\ast]$
if $c_s$ is at any point found to be not $\preceq$-minimal in $D_x$
by Lemma~\ref{lem_trivial_implies_prefix}.
\end{algo}

\begin{remark}
\label{rem_compute_minimal_elements}
A superset of $C_x$ whose cardinality is bounded by the number of atoms of
$M$ can be computed using Algorithm \ref{alg_compute_minimal} with
\texttt{f = true}, by letting $s$ range over all atoms of $M$.
Obvious short-cuts, similar to the ones described in \cite{braid_conjugacy},
can be used to increase the efficiency of this process.
\end{remark}

\section{Practical Comparisons}
\label{sect_test}
In this section, we present empirical results for braid groups $B_n$
given by the presentation (\ref{eqn_artin_pres}) from
Section \ref{sect_garside_groups}.

For several values of $n$ and $r$, we consider a set of elements $x \in B_n$
with $\slen(x) = r$, chosen at random, and compute for each such $x$ its
super summit set $S_x$ and its ultra summit set $U_x$.  Let $t_S$ and $t_U$
be the times spent on computing $S_x$ and $U_x$, respectively, and let
$n_U$ be the number of trajectories under cycling of which $U_x$ consists.
We compare the average and maximal values of $|S_x|$, $|U_x|$, $t_S$,
$t_U$ and $n_U$.
(See Tables \ref{table_uss_vs_sss_small_n} and
\ref{table_uss_vs_sss_large_n}.)

Random elements for these tests were obtained as follows.
We choose independent random simple elements $A_1,A_2,\dots$ until
$\len(A_1\cdots A_m) = r$, choose a random integer $k\in\{0,1\}$ and
compute $x = \delta^k\cdot A_1\cdots A_m$. We repeat this process until
$x$ satisfies $\slen(x) = r$.
(See Remark \ref{rem_justify_random_generator}.)
Note that $\delta^2$ is central in $B_n$,
whence there is a natural isomorphism of the graphs $\Gamma_x$ and
$\Gamma_{\delta^{2m} x}$ for arbitrary $m$. Our choice of $k$ thus
is no restriction.

In a second series of tests we consider for several values of $n$ and $r$
a set of elements $x = \delta^k\cdot A_1\cdots A_r \in B_n$ obtained by
choosing a random integer $k\in\{0,1\}$ and independent random
simple elements $A_1,\dots,A_r$.  We compare the average values of
$\len(x)$ and $\slen(x)$, as well as the percentages $\epsilon_S$ and
$\epsilon_U$ of elements $x$ satisfying $x\in S_x$ and $x\in U_x$,
respectively. (See Table \ref{table_density}.)

All computations were performed on a Linux PC with a 2.4\,GHz Pentium
4 CPU, 533\,MHz system bus and  512\,MB of RAM using the author's
implementation in C, which is part of the computational algebra system
\Magma{}~\cite{magma}.

\begin{table}
 \caption{Average / maximal values for $|U_x|$, $|S_x|$, $t_U$, $t_S$
 and $n_U$ (see text); times are given in ms, unless stated otherwise.
 Where no values of $|S_x|$ and $t_S$ are given, computing super
 summit sets exceeded the available memory of 512\,MB.}
 \label{table_uss_vs_sss_small_n}
 \medskip
 \begin{tabular}{@{$\:$}c@{$\:$}||@{$\:$}c@{$\:$}|@{$\:$}c@{$\:$}|@{$\:$}c@{$\:$}|@{$\:$}c@{$\:$}|@{$\:$}c@{$\:$}|@{$\:$} c@{$\:$}}
  $n$     & \multicolumn{6}{c}{3} \\
  \hline
  $r$     & 2 & 5 & 10 & 20 & 100 & 1000 \\
  \hline
  $|U_x|$ & 3.1 / 4 & 9.8 / 10 & 20 / 20 & 40 / 40 &
            200 / 200 & 2000 / 2000  \\
  $|S_x|$ & 3.1 / 4 & 9.8 / 10 & 20 / 20 & 40 / 40 &
            200 / 200 & 2000 / 2000  \\
  $t_U$   & 0.1 / 10 & 0.2 / 10 & 0.4 / 11 & 1.1 / 11 &
            22 / 31 & 4.1\,s / 5.4\,s  \\
  $t_S$   & 0.1 / 9 & 0.3 / 10 & 1.0 / 11 & 3.4 / 11 &
            79 / 90 & 15\,s / 19\,s  \\
  $n_U$   & 1.2 / 2 & 1.5 / 2 & 1.5 / 2 & 1.5 / 2 &
            1.4 / 2 & 1.6 / 2  \\
\multicolumn{7}{c}{}\\ 
  $n$     & \multicolumn{6}{c}{4}  \\
  \hline
  $r$     & 2 & 5 & 10 & 20 & 100 & 1000 \\
  \hline
  $|U_x|$ & 5.6 / 10 & 12 / 50 & 20 / 40 & 40 / 40 &
            200 / 200 & 2000 / 2000  \\
  $|S_x|$ & 11 / 24 & 47 / 128 & 100 / 464 & 190 / 660 &
            920 / 1704 & 9000 / 1.0e4 \\
  $t_U$   & 0.2 / 11 & 0.5 / 11 & 0.7 / 11 & 1.8 / 11 &
            45 / 81 & 7.8\,s / 13.5\,s \\
  $t_S$   & 0.4 / 11 & 2.6 / 11 & 9.2 / 51 & 29 / 121 &
            650 / 1250 & 210\,s / 272\,s \\
  $n_U$   & 1.6 / 3 & 1.7 / 10 & 1.5 / 8 & 1.5 / 2 &
            1.5 / 2 & 1.6 / 2  \\
\multicolumn{7}{c}{}\\ 
  $n$     & \multicolumn{6}{c}{6} \\
  \hline
  $r$     & 2 & 5 & 10 & 20 & 100 & 1000 \\
  \hline
  $|U_x|$ & 15 / 72 & 17 / 1440 & 21 / 60 & 40 / 40 &
            200 / 200 & 2000 / 2000  \\
  $|S_x|$ & 270 / 1004 & 3800 / 8.3e4 & 1.1e4 / 2.9e5 & --- &
            --- & ---  \\
  $t_U$   & 1.3 / 11 & 1.9 / 151 & 1.6 / 30 & 3.1 /  20 &
            53 / 90 & 5.2\,s / 12\,s \\
  $t_S$   & 18 / 71 & 600 / 15\,s & 24\,s / 672\,s & --- &
          --- & ---  \\
  $n_U$   & 3.1 / 18 & 2.6 / 262 & 1.5 / 4 & 1.5 / 2 &
            1.4 / 2 & 1.6 / 2  \\
\multicolumn{7}{c}{}\\ 
  $n$     & \multicolumn{6}{c}{8} \\
  \hline
  $r$     & 2 & 5 & 10 & 20 & 100 & 1000 \\
  \hline
  $|U_x|$ & 43 / 448 & 14 / 188 & 21 / 56 & 40 / 40 &
            200 / 200 & 2000 / 2000  \\
  $|S_x|$ & 1.3e4 / 7.3e4  & --- & --- & --- & --- & ---  \\
  $t_U$   & 4.9 / 59 & 2.5 / 80 & 1.9 / 40 & 4.7 / 11 &
            67 / 150 & 7.7\,s / 17\,s \\
  $t_S$   & 27\,s / 165\,s & --- & --- & --- & --- & ---  \\
  $n_U$   & 6.9 / 64 & 2.7 / 94 & 1.5 / 2 & 1.5 / 2 &
            1.5 / 2 & 1.4 / 2
 \end{tabular}
 \medskip
\end{table}

\begin{table}
 \caption{Average / maximal values for $|U_x|$, $|S_x|$, $t_U$, $t_S$
 and $n_U$ (see text); times are given in ms, unless stated otherwise.
 In all cases, computing super summit sets exceeded the available
 memory of 512\,MB.}
 \label{table_uss_vs_sss_large_n}
 \medskip
 \begin{tabular}{c||c|c|c|c|c|c}
  $n$     & \multicolumn{6}{c}{10} \\
  \hline
  $r$     & 2 & 5 & 10 & 20 & 100 & 1000 \\
  \hline
  $|U_x|$ & 63 / 1408 & 15 / 54 & 21 / 40 & 40 / 78 &
            200 / 200 & 2000 / 2000  \\
  $t_U$   & 12 / 290 & 3.3 / 21 & 4.2 / 40 & 6.3 / 90 &
            100 / 190 & 16\,s / 32\,s \\
  $n_U$   & 11 / 104 & 2.0 / 8 & 1.5 / 4 & 1.6 / 2 &
            1.5 / 2 & 1.5 / 2  \\
\multicolumn{7}{c}{}\\ 
  $n$     & \multicolumn{6}{c}{20} \\
  \hline
  $r$     & 2 & 5 & 10 & 20 & 100 & 1000 \\
  \hline
  $|U_x|$ & 30 / 280 & 12 / 20 & 20 / 40 & 40 / 40 &
            200 / 200 & 2000 / 2000  \\
  $t_U$   & 10 / 151 & 3.4 / 11 & 4.7 / 11 & 9.7 / 21 &
            100 / 221 & 19\,s / 46\,s \\
  $n_U$   & 7.7 / 70 & 1.9 / 4 & 1.5 / 4 & 1.5 / 2 &
            1.6 / 2 & 1.5 / 2  \\
\multicolumn{7}{c}{}\\ 
  $n$     & \multicolumn{6}{c}{50} \\
  \hline
  $r$     & 2 & 5 & 10 & 20 & 100 & 1000 \\
  \hline
  $|U_x|$ & 7.0 / 64 & 10 / 20 & 20 / 20 & 40 / 40 &
            200 / 200 & 2000 / 2000  \\
  $t_U$   & 7.8 / 50 & 8.4 / 21 & 12 / 21 & 18 / 30 &
            130 / 241 & 21\,s / 48\,s \\
  $n_U$   & 2.3 / 16 & 1.6 / 4 & 1.5 / 2 & 1.5 / 2 &
            1.5 / 2 & 1.6 / 2  \\
\multicolumn{7}{c}{}\\ 
  $n$     & \multicolumn{6}{c}{100} \\
  \hline
  $r$     & 2 & 5 & 10 & 20 & 100 & 1000 \\
  \hline
  $|U_x|$ & 5.2 / 32 & 10 / 10 & 20 / 20 & 40 / 40 &
            200 / 200 & 2000 / 2000  \\
  $t_U$   & 20 / 101 & 27 / 50 & 36 / 61 & 49 / 69 &
            210 / 370 & 23\,s / 32\,s  \\
  $n_U$   & 1.7 / 8 & 1.4 / 2 & 1.5 / 2 & 1.6 / 2 &
            1.5 / 2 & 1.5 / 2
 \end{tabular}
 \medskip
\end{table}

\begin{table}
 \caption{Average values of $\len(x)$ and $\slen(x)$ and percentages
 $\epsilon_S$ and $\epsilon_U$ of elements satisfying $x\in S_x$ and
 $x\in U_x$, respectively. (See text.)}
 \label{table_density}
 \medskip
 \begin{tabular}{@{$\:$}c||c|c|c|c|c|c||c|c|c|c|c|c@{$\:$}}
  $n$            & \multicolumn{6}{c||}{3}
                   & \multicolumn{6}{c}{4} \\
  \hline
  $r$            & 2 & 5 & 10 & 20 & 100 & 1000
                   & 2 & 5 & 10 & 20 & 100 & 1000 \\
  \hline
  $\len(x)$      & 1.0 & 1.8 & 2.7 & 4.7 & 19 & 170
                   & 1.4 & 2.7 & 4.5 & 7.8 & 34 & 330 \\
  $\slen(x)$    & 0.8 & 1.4 & 2.1 & 3.7 & 17 & 170
                   & 1.2 & 2.1 & 3.6 & 6.6 & 33 & 330 \\
  $\epsilon_S$   & 89 & 72 & 64 & 56 & 52 & 51
                   & 77 & 53 & 41 & 36 & 32 & 32 \\
  $\epsilon_U$   & 89 & 72 & 64 & 56 & 52 & 51
                   & 72 & 40 & 22 & 11 & 8.7 & 8.0 \\
\multicolumn{7}{c}{}\\ 
  $n$            & \multicolumn{6}{c||}{6}
                   & \multicolumn{6}{c}{10} \\
  \hline
  $r$            & 2 & 5 & 10 & 20 & 100 & 1000
                   & 2 & 5 & 10 & 20 & 100 & 1000 \\
  \hline
  $\len(x)$      & 1.9 & 3.8 & 6.7 & 12 & 58 & 570
                   & 2.0 & 4.8 & 9.0 & 17 & 85 & 840 \\
  $\slen(x)$    & 1.6 & 3.1 & 5.6 & 11 & 57 & 570
                   & 2.0 & 4.3 & 8.4 & 17 & 84 & 840 \\
  $\epsilon_S$   & 77 & 42 & 33 & 32 & 32 & 31
                   & 96 & 63 & 55 & 51 & 54 & 53 \\
  $\epsilon_U$   & 30 & 4.0 & 1.1 & 0.9 & 1.0 & 1.4
                   & 1.4 & 0.3 & 0.0 & 0.0 & 0.1 & 0.0 \\
\multicolumn{7}{c}{}\\ 
  $n$            & \multicolumn{6}{c||}{15}
                   & \multicolumn{6}{c}{30, 50, 75, 100} \\
  \hline
  $r$            & 2 & 5 & 10 & 20 & 100 & 1000
                   & 2 & 5 & 10 & 20 & 100 & 1000 \\
  \hline
  $\len(x)$      & 2.0 & 5.0 & 9.9 & 20 & 98 & 980
                   & 2.0 & 5.0 & 10 & 20 & 100 & 1000 \\
  $\slen(x)$    & 2.0 & 4.9 & 9.8 & 20 & 98 & 980
                   & 2.0 & 5.0 & 10 & 20 & 100 & 1000 \\
  $\epsilon_S$   & 100 & 94 & 89 & 87 & 88 & 87
                   & 100 & 100 & 100 & 100 & 100 & 100 \\
  $\epsilon_U$   & 0.0 & 0.0 & 0.0 & 0.1 & 0.0 & 0.0
                   & 0.0 & 0.0 & 0.0 & 0.0 & 0.0 & 0.0
 \end{tabular}
 \medskip
\end{table}

\subsection{Results}
\label{sect_test_result}

The main results of the tests can be summarised as follows.
\begin{enumerate}
\item The average size of $S_x$ grows very fast with increasing values of
 $n$.  $S_x$ is in general not computable on typical current computers for
 $n\ge 10$ or $n>5$, $r>15$, due to extreme memory requirements.

\item With the exception of very small values of $r$ ($r=2,5$),
 the average size of $U_x$ is of the order of $2r$, in particular almost
 independent of $n$, for the case of presentation (\ref{eqn_artin_pres})
 from Section \ref{sect_garside_groups}.  Similar tests for presentation
 (\ref{eqn_bkl_pres}) yield an average size of the order of $nr$ for not
 too small values of $r$.

 There are, however, elements whose ultra summit sets are much larger than
 the average values.  With growing values of $n$ and $r$, these exceptions
 seem to get rarer, so in some sense the situation then becomes easier.

 In the tests, $U_x$ remained sufficiently small to be computed easily
 over the entire parameter range.
 
\item The average number of connected components (trajectories) of $U_x$
 is approximately 1.5 for larger values of $r$.  Note that this implies
 that computing conjugating elements by Algorithm \ref{alg_conjugacy_search}
 is very efficient.

 Another consequence of this is that even in the case $n = 3$ where
 $U_x = S_x$,
 computing $U_x$ is much faster than computing $S_x$ for large values
 of $r$, since the decomposition of $U_x$ into trajectories is used
 efficiently (Theorem \ref{thm_minimal_simple_elements}).

\item A random element of the form $\delta^k\cdot A_1\cdots A_r$ with
 independent random simple elements $A_1,\dots,A_r$ is surprisingly likely
 to be a super summit element, that is, satisfy $x\in S_x$. In the tests
 for $n>20$, the probability for this is indistinguishable from 1 and
 the elements moreover satisfy $\slen(x) = r$.

 Random elements as above which are ultra summit elements, on the other
 hand, are very rare for $n>5$ and were not encountered at all in the
 tests for $n>20$.

 This suggests that, with the exception of braid groups on very few
 strings, the ultra summit set of an element in general is a very small
 subset of the super summit set.
\end{enumerate}

\begin{remark}
\label{rem_justify_random_generator}
Other methods of constructing random elements may
produce different distributions on the set of all elements $x \in B_n$
with $\slen(x) = r$ and $x\in S_x$.  However, the most natural (and
computationally most efficient) way of producing random elements
seems to be computing the normal forms of longer sequences of
independent random simple elements and hence this method was used
in the tests.

Note, moreover, that at least for larger values of $n$, according to
our results a product $x$ of a random power of $\delta$ and $r$
independent random simple elements is extremely likely to satisfy both
$x\in S_x$ and $\len(x)=\slen(x)=r$.  In this sense, the distribution
of random super summit elements with given canonical length produced
by the method used in our tests is very natural.

According to tests with other methods of generating random elements,
our main results as formulated in Section \ref{sect_test_result} in
any case do not seem to depend on the details of random element
generation.
\end{remark}

\section{Conclusions}
\label{sect_conclusions}

We defined a new invariant of conjugacy classes in Garside groups, the
ultra summit set, using the digraph structure of the well-known super
summit set induced by the cycling operation and established that it
satisfies \lq\lq{}convexity\rq\rq{} properties analogous to the ones
holding for super summit sets.  Ultra summit sets seem to be rather natural
objects and may be useful for further theoretical analysis of Garside groups.

Apart from their theoretical significance, our results allow efficient
computation of ultra summit sets, providing a practical solution to the
conjugacy decision and search problems in Garside groups.

Our tests for Artin's presentation of $B_n$ show that, in particular
for larger braid index $n$, super summit elements are extremely common
and super summit sets hence are much too large to be of computational use.
Ultra summit elements, on the other hand, seem to be extremely rare and
ultra summit sets can be computed easily even for large values of braid
index and canonical length.  We demonstrate that, using ultra summit sets,
the conjugacy decision and search problems can be solved in very little
time on current computers for elements of canonical length 1000 in $B_{100}$.

Hence from both a theoretical and a computational point of view, the notion of
ultra summit sets appears to be a significant advance in the study of the
conjugacy problems in Garside groups.

\section*{Acknowledgment}
I am grateful to Prof.~John Cannon for his support for this research project.
I would like to thank Prof.~Patrick Dehornoy and Prof.~Hugh Morton for
valuable comments on an earlier version of this paper.

\bibliographystyle{elsart-num}
\bibliography{literatur}

\end{document}